\documentclass[letterpaper,12pt]{article}

\usepackage{graphicx}
\usepackage{amssymb,amsmath}
\usepackage{subfigure}
\usepackage{amsthm}
\usepackage{hyperref}
\usepackage{color}
\usepackage{multirow,txfonts,verbatim}
\usepackage{mathrsfs}
\usepackage{booktabs}
\usepackage{bm}
\usepackage{tikz}
\usetikzlibrary{decorations.pathreplacing,calligraphy}

\usepackage{algorithm}
\usepackage{algpseudocode}
\usepackage[square,comma,sort&compress,numbers]{natbib}

\usepackage{float}

\newtheoremstyle{plainNoItalics}{}{}{\normalfont}{}{\bfseries}{.}{ }{}
\theoremstyle{plain}

\theoremstyle{plainNoItalics}

\newtheorem{example}{Example}[section]
\theoremstyle{definition}
\theoremstyle{definition}
\theoremstyle{definition}
\theoremstyle{definition}
\theoremstyle{definition}
\theoremstyle{definition}
\numberwithin{equation}{section}
\numberwithin{theorem}{section}
\numberwithin{lemma}{section}
\numberwithin{figure}{section}
\numberwithin{example}{section}
\numberwithin{rem}{section}

\usepackage{lineno}
\usepackage{xcolor}
\usepackage{natbib}

\newcommand{\bU}{{\bm{U}}}
\newcommand{\bF}{{\bm{F}}}

\newcommand{\bxi}{\bm{\xi}}
\newcommand{\bd}{{\bm{d}}}
\newcommand{\bE}{{\bm{E}}}

\newcommand{\bW}{{\bm{W}}}

\graphicspath{{figures/}}

\title{A multi-fidelity machine learning based semi-Lagrangian finite volume scheme for linear transport equations and the nonlinear Vlasov-Poisson system}
\author{
Yongsheng Chen
\thanks{School of Mathematical Sciences, Zhejiang University, Hangzhou, 310027, China. {\tt 22035024@zju.edu.cn}}
\and
Wei Guo
\thanks{Corresponding author. Department of Mathematics and Statistics, Texas Tech University, Lubbock, TX, 70409, USA. 
{\tt weimath.guo@ttu.edu}. }
\and
Xinghui Zhong
\thanks{School of Mathematical Sciences, Zhejiang University, Hangzhou, 310027, China. {\tt zhongxh@zju.edu.cn}}
}

\begin{document}

\maketitle
\begin{abstract}

Machine-learning (ML) based discretization has been developed to simulate complex partial  differential equations (PDEs) with tremendous success   across various fields. These learned PDE solvers can effectively resolve the underlying solution structures of interest and achieve a level of accuracy which often requires an order-of-magnitude finer grid for a  conventional numerical method using polynomial-based approximations. In a previous work in \cite{chen_learned_2023}, we introduced a learned 
finite volume discretization that further incorporates the semi-Lagrangian (SL) mechanism, enabling larger CFL numbers for stability. However, the efficiency and effectiveness of such methodology heavily rely on the availability of abundant high-resolution training data, which can be prohibitively expensive to obtain. To address this challenge, in this paper, we propose a novel multi-fidelity ML-based SL method for transport equations. This method leverages a combination of a small amount of high-fidelity data and sufficient but cheaper low-fidelity data. The approach is designed based on a composite convolutional neural network architecture that explore the inherent correlation between high-fidelity and low-fidelity data.   The proposed method demonstrates the capability to achieve a reasonable level of accuracy, particularly in scenarios where a single-fidelity model fails to generalize effectively. We further extend the method to the nonlinear Vlasov-Poisson system by employing high order Runge-Kutta exponential integrators.  
A collection of numerical tests are provided to validate the efficiency and accuracy of the proposed method.

\end{abstract}
\section{Introduction}

The rapid development of machine learning (ML) has opened new research avenues for approximating complex partial differential equations (PDEs), and many successful ML-based PDE solvers are designed by leveraging the  expressive power of neural networks (NNs) and advancements of automatic differentiation technology \cite{Ilya_fix_2017}.
 Among these developments, ML-based discretization has received substantial research attention. This approach combines the strengths of ML techniques and conventional numerical methods, leading to  remarkable success across diverse scientific and engineering applications \cite{obiols2020cfdnet,sirignano2020dpm,pathak2020using,um2020solver,tompson2017accelerating}. Such 
equation-specific ML-based discretization replaces the polynomial-based approximation with NNs, yielding more flexible and accurate representations of the underlying PDEs. With high-fidelity data, optimal numerical discretization can be learned through the training process of NNs. This enables such ML-based discretization to achieve accurate and stable numerical solutions even with much coarser grid resolutions \cite{bar-sinai_learning_2019,zhuang2021learned,kochkov_machine_2021}
, reducing computational cost compared to traditional numerical methods that require finer grid for comparable accuracy. Other approaches of ML-based PDE solvers include the physics informed neural networks (PINNs) which utilize the physics-informed loss  
\cite{xiang_self-adaptive_2022,yu_gradient-enhanced_2022,lu_physics-informed_2021,raissi_physics_2017,raissi_physics_2017-1}, 
neural operators
\cite{lu_deeponet_2021,li_fourier_2021,kissas_learning_2022,li_neural_2020,bhattacharya_model_2021,trask_gmls-nets_2019}
, and autoregressive methods 
\cite{bar-sinai_learning_2019,greenfeld2019learning,hsieh2019learning,brandstetter2022message,chen_cell-average_2022}
, etc.. Compared to these methods, ML-based discretization offers an additional advantage: the ability to enforce inherent physical constraints, such as conservation of mass, momentum, and energy of the system, at the discrete level. The integration of such inductive biases plays a pivotal role in enhancing the generalization capabilities of the NN 
\cite{bar-sinai_learning_2019,zhuang2021learned,kochkov_machine_2021}.

In our recent work \cite{chen_learned_2023}, we proposed an autoregressive ML finite volume (FV) 
 method under the ML-based discretization framework. This method couples with the semi-Lagrangian (SL) mechanism for solving linear transport equations. By incorporating a specific inductive bias into the NN, the method aims to learn the SL discretization from the data, avoiding the need for costly upstream cells tracking. The proposed method inherits all the advantages of the ML-based discretization approach. Additionally, it allows for a larger Courant-Friedrichs-Lewy (CFL) number for stability compared to the Eulerian method-of-line approach, thereby enhancing the efficiency. However, the success of this method heavily relies on the availability of sufficient high-fidelity data, which can be prohibitively expensive to generate, especially for complex PDE simulations. This is known as a primary challenge for ML-based methods. While high-fidelity data is ideal for
 training the NN, but it is often difficult to acquire and scarce in quantity. On the other hand, the low-fidelity data is easily obtainable but is insufficient for training ML models with satisfactory accuracy.
 Therefore, using a single-fidelity model is prone to generalization failure.  To address this limitation, it becomes necessary to leverage data with multiple fidelity levels for multi-fidelity modeling. This approach can enhance both accuracy and generalization capability of the model \cite{fernandez2016review,peherstorfer_survey_2018}.

In the literature,  multi-fidelity modeling with NNs can generally be categorized into three approaches. The first approach focuses on learning the correlation between the low-fidelity and high-fidelity data. For instance, it has been employed to approximate the linear or nonlinear correlations between low-fidelity and high-fidelity solutions for PINNs \cite{liu_multi-fidelity_2019,meng_composite_2020,ramezankhani_data-driven_2022} 
and for operator learning using the DeepONet approach 
\cite{lu_multifidelity_2022,de_bi-fidelity_2023,howard_multifidelity_2022}
. Multi-fidelity Bayesian NNs 
\cite{meng_multi-fidelity_2021}
, which integrate the Bayesian framework and PINNs, provide an example of capturing the cross-correlation with uncertainty quantification between the low- and high-fidelity data. Recently, the authors in 
\cite{chen_multi-fidelity_2022}
proposed a convolutional NN (CNN) architecture
to exploit the correlation among the multi-fidelity data. 
More recently, the authors in \cite{chen2023feature} 
proposed a multi-fidelity PINN, where the low- and high-fidelity solutions are projected onto the same feature space using an encoder, and their projections are adjacent by constraining their relative distance with the NN.

The second approach in multi-fidelity learning involves utilizing transfer learning \cite{weiss_survey_2016}. 
In this approach, the NN is initially trained with low-fidelity information and subsequently fine-tuned with high-fidelity information, aiming to improve the overall accuracy of the model. Several related works on PDE solving tasks with transfer learning can be found in \cite{aliakbari_predicting_2022,ashouri_transfer_2022,chakraborty_transfer_2021,chen_transfer_2021,de_transfer_2020,de_neural_2022,jiang_use_2023,song_transfer_2022,motamed_multi-fidelity_2020,li_line_2022}. 
The third approach involves learning low-fidelity and high-fidelity solutions using the same NN. This approach can be achieved by using the low-fidelity solution as an intermediate value of the network 
\cite{guo_multi-fidelity_2022}
or by imposing constraints on the solution with any available high-fidelity data \cite{basir_physics_2022}.


In this paper, we propose a novel multi-fidelity ML-based method for the ML-based SL FV method solving transport equations, with the focus on scenarios where there is an abundance of low-fidelity data and limited high-fidelity data. 
Our method belongs to the first category and aims to take advantage of the accuracy
of the high-fidelity data and the accessibility of the low-fidelity data simultaneously. To achieve this, we introduce a composite NN architecture 
that consists of two sub-networks: the low-fidelity network and the high-fidelity network. The low-fidelity network shares the same structure as the ML-based SL FV scheme presented in 
\cite{chen_learned_2023},
and it is intended  to predict low-fidelity solutions. The high-fidelity network also has a similar structure but incorporate the output of the low-fidelity network as an additional input to approximate the correlation between the low-fidelity and high-fidelity solutions. By integrating the two networks, our method achieves improved stability and accuracy compared to using networks trained solely on low-fidelity or high-fidelity data. It also exhibits satisfactory generalization capabilities. Moreover, our method inherits the advantages of the ML-based SL FV scheme 
\cite{chen_learned_2023}
, such as mass conservation, translational equivariance, avoiding the need to tracking upstream cells,  the ability to allow for large CFL numbers, and attaining an accuracy level that exceeds that of traditional numerical algorithms with the same mesh resolution.

In addition to our proposed multi-fidelity ML-based SL FV method for solving linear transport equations, we also propose to extend this method for simulating the nonlinear Vlasov-Poisson (VP) system. The VP system presents an additional challenge for accurate tracking of its characteristic equations due to the inherent nonlinearity.  To overcome this challenge, we combine the high-order Runge-Kutta (RK) exponential integrators (RKEI),  introduced in \cite{celledoni2003commutator,cai_high_2021}.
By employing the RKEI, the VP system can be decomposed into a sequence of linearized transport equations with frozen coefficients \cite{cai_high_2021,zheng2022fourth}. 
This decomposition allows us to apply the multi-fidelity model, resulting in a nonlinear data-driven multi-fidelity SL FV scheme for the VP system. 

The rest of the paper is organized as follows. In Section \ref{sec:algorithm}, we review the ML-based SL FV scheme proposed in \cite{chen_learned_2023}. 
Section \ref{sec:multi} is devoted to presenting our multi-fidelity data-driven SL FV scheme. In particular, we lay out the details for linear transport equations in Section \ref{sec:scheme:transport} and discuss its extension to the nonlinear VP system by coupling with the RKEI method in Section \ref{sec:vp}. 
In Section \ref{sec:num}, numerical results are provided to demonstrate the performance of our multi-fidelity method. The conclusion and future work are discussed in Section \ref{sec:con}.

\section{ML-based SL FV scheme}\label{sec:algorithm}
In this section, we provide a brief overview of the ML-based SL method for solving linear transport equations developed in \cite{chen_learned_2023}. 
The methodology will  serve as a building block of the proposed multi-fidelity solver discussed in Section \ref{sec:multi}.   

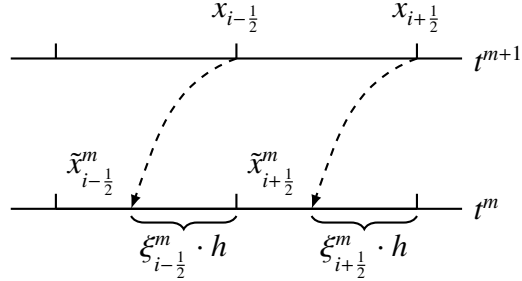
\begin{figure}[!htbp]
	\centering
		\begin{tikzpicture}[x=1cm,y=1cm]
		\begin{scope}[thick]
		\draw[black]  (-3,0) -- (3,0)
		node[right]{$t^{m}$};
		\draw[black] (-3,2) -- (3,2)
		node[right]{$t^{m+1}$};
		
		
		\draw[black] (-2.4,0) -- (0,0);
		\draw[black] (0,0) -- (2.4,0);
		\draw[black] (0,2) -- (0,2.2)node[above] {$x_{i-\frac12}$};
            \draw[black] (0,0) -- (0,0.2);
		\draw[black] (2.4,2) -- (2.4,2.2)node[above] {$x_{i+\frac12}$};
            \draw[black] (2.4,0) -- (2.4,0.2);
		\draw[black] (-2.4,2) -- (-2.4,2.2);
            \draw[black] (-2.4,0) -- (-2.4,0.2);            
		\draw[-latex,dashed](0,2) to[out=200,in=70] (-1.4,0) node[above left=1pt] {$\tilde{x}^m_{i-\frac12}$};  
		\draw[-latex,dashed](2.4,2) to[out=200,in=70] ( 1.,0) node[above left=1pt] {$\tilde{x}^m_{i+\frac12}$};
  \draw [decorate,decoration={brace,amplitude=5pt,mirror,raise=0.5ex}]
  (-1.4,0) -- (0,0)node[midway,yshift=-1.5em] {$\xi^m_{i-\frac12}\cdot h$};
    \draw [decorate,decoration={brace,amplitude=5pt,mirror,raise=0.5ex}]
  (1,0) -- (2.4,0)node[midway,yshift=-1.5em] {$\xi^m_{i+\frac12}\cdot h$};
		
		\end{scope}
		\end{tikzpicture}

	\caption{Schematic illustration of the 1D SL FV scheme.\label{fig:slfv}}	
\end{figure}

Consider the following one-dimensional (1D) transport equation
\begin{equation}
    \label{eq:transport1d}
    u_t + (a(x,t)u)_x = 0,\quad x\in\Omega,
\end{equation}
where $a(x,t)$ denotes the velocity field. The domain $\Omega$ is uniformly partitioned into $N$ cells, where each cell is represented by  
$I_i=[x_{i-\frac12},x_{i+\frac12}]$, and the mesh size is denoted by $h = x_{i+\frac12} - x_{i-\frac12}$. 
In an SL FV scheme, the solution is represented by the cell-averaged values $U_i^m$ at time $t^m$ in the cell $I_i$. 
The scheme then advances $\{U_i^m\}$ using characteristic information. Specifically, the scheme traces the characteristics backwards in time
from time step $t^{m+1}$ to $t^m$ by
\begin{equation}
\label{eq:characteristic}
\begin{cases}
\displaystyle\frac{dx(t)}{dt} &= a(x(t),t),\\[1mm]
x(t^{m+1})  &= x_{i+\frac12},\;\;i=1,\ldots,N,
\end{cases}
\end{equation}
to obtain the endpoints of the upstream cell, denoted by $\tilde{x}^m_{i+\frac12} = x(t^m)$, as shown in Figure \ref{fig:slfv}.  The advantage of integrating characteristic information in the algorithm design   is that it allows for larger CFL number and achieves higher accuracy and efficiency.
Then $\{U_i^{m+1}\}$ can be updated based on the fact that the exact solution $u(x,t)$ satisfies
\begin{align*}
    \int_{I_i} u\left(x,t^{m+1}\right) dx = \int_{\tilde{x}^m_{i-\frac12}}^{\tilde{x}^m_{i+\frac12}} u\left(x,t^{m}\right)dx,\;\;i=1,\ldots,N,
\end{align*} 
and the right hand side can be approximated by integrating 
 the polynomials  in $[\tilde{x}^m_{i-\frac12},\tilde{x}^m_{i+\frac12}]$ reconstructed using the standard FV framework. With the normalized shift defined by
$$
\xi^m_{i+\frac12} = \frac{\tilde{x}^m_{i+\frac12} - x_{i+\frac12}}{h},\quad  i=1,\ldots,N,
$$
the SL FV schemes  can be rewritten as
the following formulation
\begin{equation}
\label{eq:sl_re}
U_i^{m+1} = \sum_{\ell\in\mathcal{S}^m_i} d^m_{i,\ell} U_{\ell}^m,
\end{equation}
where $\mathcal{S}^m_i$ denotes the stencil used for updating  $U_i^{m+1}$, and $\{d_{i,\ell}^m\}$ are the associated coefficients fully determined by  $\{U_i^m\}$ and $\{\xi_i^m\}$, and more details can be found in \cite{qiu_conservative_2011}.
Moreover, the scheme \eqref{eq:sl_re} is provably mass conservative when
\begin{equation}
\sum_{i}d_{i,\ell}^m = 1,\quad\forall \ell.\label{eq:mass}
\end{equation}

\begin{figure}[!htbp]
\centering
    \centerline{\includegraphics[width=0.95\textwidth]{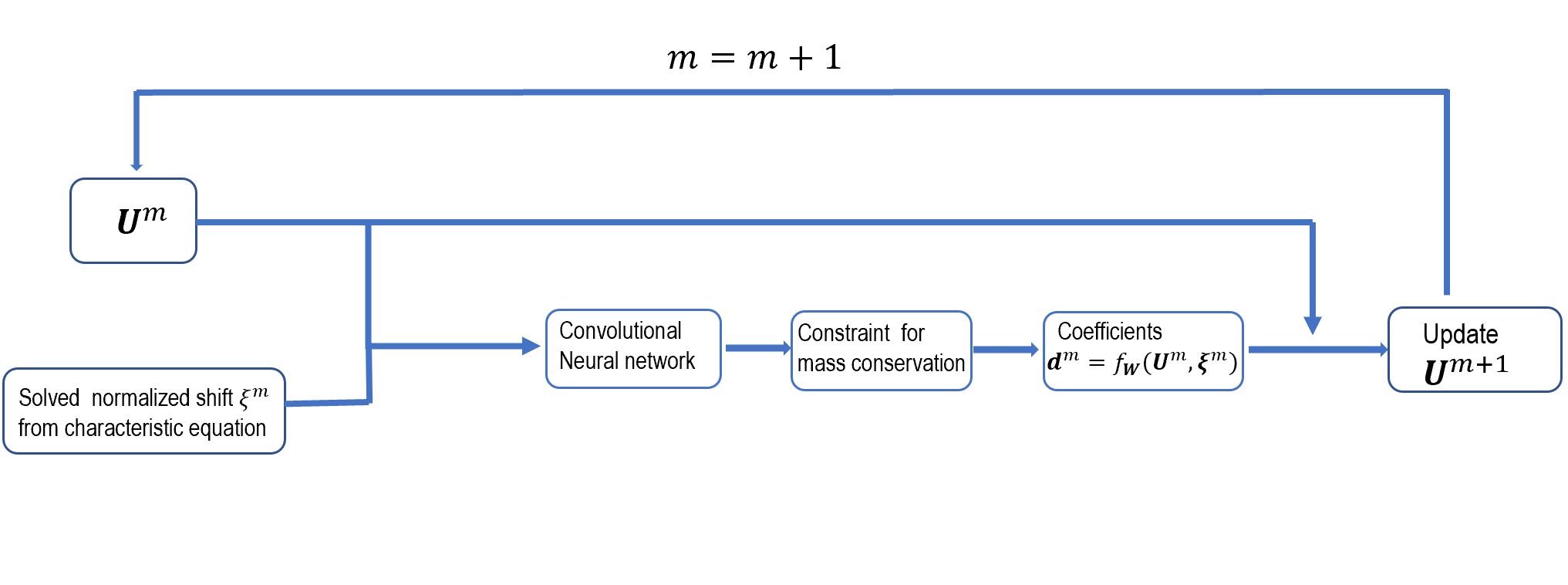}}
	\caption{Illustration of the proposed ML-based SL FV method in \cite{chen_learned_2023}
.}\label{fig:network_structure}
\end{figure}

It is worth emphasizing that the computation of $\{d_{i,\ell}^m\}$ requires expensive tracking of the geometries for upstream cells. This procedure not only demands sophisticated implementation but also significantly dominates the overall computational cost, especially in high dimensions \cite{lauritzen2010conservative,erath2013mass}. 
To address this challenge, we proposed an ML-based SL FV scheme in 
\cite{chen_learned_2023},
for which the discretization coefficients $\{d^m_{i,\ell}\}$ are learned from the data, drawing inspiration from related works 
\cite{bar-sinai_learning_2019,zhuang2021learned,kochkov_machine_2021}. 
In particular, the proposed scheme incorporates the normalized shifts $\{\xi_i^m\}$ as an inductive bias and 
uses the NN $f_\bW$ to infer the discretization coefficients:
 \begin{equation}
 \label{eq:fw}
    \bd^m = f_\bW(\bU^m,\bxi^m),
\end{equation} 
where $\bU^m$, $\bxi^m$, and $\bd^m$ denote the collection of $\{U^m_i\}$, $\{\xi^m_{i-\frac12}\}$, and $\{d^m_{i,\ell}\}$, respectively. The NN $f_\bW$ takes $\bU^m$ and  $\bxi^m$ as two-channel inputs and is constructed as a stack of convolutional layers with trainable parameters $\bW$ and nonlinear activation functions, such as ReLU. Additionally, a constraint layer is incorporated to enforce \eqref{eq:mass}, ensuring exact mass conservation. Once $\bm{d}^m$ is obtained, the solution $\bU^{m+1}$ is updated with \eqref{eq:sl_re}. The scheme is schematically illustrated in Figure \ref{fig:network_structure}.

Unlike the traditional SL FV schemes, this ML-based scheme unitizes a set of fixed centered stencils. For example, we can choose the 5-cell stencil $\mathcal{S}^m_{i} = \{i-2,i-1,i,i+1,i+2\}$ to update the solution $U_i^{m+1}$. In this way, 
this approach offers the benefits of simplifying algorithm development and conveniently satisfying the mass conservation constraint, though it compromises the unconditional stability of traditional SL FV scheme, which relies on local stencils for reconstruction. Numerical evidences demonstrate that the ML-based SL FV scheme can maintain stability for CFL numbers up to 1.8. In addition, the scheme features the translational equivariance, which is highly desirable for data efficiency and improved generalization capability. Further, the method can be easily adapted for  two-dimensional (2D) linear transport problems with slight modification. Within this framework, the enforcement of physical constraints, such as the conservation of mass of the system, can be conveniently achieved at the discrete level, distinguishing it from many other ML-based PDEs solvers.

The rationale behind the ML-based PDE discretization approach lies in the observation that the solution manifold  often exhibits low dimensional structures, including recurrent patterns and  coherent structures.  By leveraging high resolution training  data, the underlying NN effectively learns an optimal discretization that can achieve significantly higher accuracy compared with a traditional polynomial-based method.  However, similar to other ML-based discretization schemes, the effectiveness of this method heavily relies on the availability of abundant high-resolution training data. In particular, to generate the training data, we sample a set of initial conditions from a prescribed distribution. From each sampled initial condition, we employ an accurate and reliable FV method, such as an SL WENO method or Eulerian WENO method, to generate solution trajectories on a high-resolution mesh. It is crucial to adequately resolve the solution structures of interest during this process. Subsequently, we obtain the training data by downsampling these solution trajectories to a mesh with reduced resolution by an order of magnitude. It is worth noting that acquiring such training data can often be prohibitively expensive, especially for complex problems such as plasma simulations. Moreover, the scarcity of high-resolution data can significantly limit the performance of the scheme and impede its ability to generalize effectively. 

\section{Multi-fidelity data-driven  SL FV scheme}\label{sec:multi}
In this section, we present a novel data-driven multi-fidelity SL FV scheme
aimed at reducing the reliance on high-resolution data for the ML-based SL FV method introduced in Section \ref{sec:algorithm}, while maintaining a reasonable level of accuracy. We begin by introducing the scheme for transport equations in Section \ref{sec:scheme:transport}, and subsequently extend the scheme to the nonlinear VP system in Section \ref{sec:vp}.
\subsection{Transport equations}
\label{sec:scheme:transport}
In this subsection, we lay out the details of the multi-fidelity approach for linear transport equation with the goal of enhancing the accuracy and efficiency of the ML-based SL FV method while minimizing the reliance on high-resolution data.

The training dataset for our multi-fidelity approach is composed of both high-fidelity and low-fidelity data. The high-fidelity data, denoted by $
\{\bU_H\}$, comprises coarsened high-resolution numerical solution trajectories, which are accurate but computationally expensive and scarce. On the other hand, the low-fidelity data, denoted by $\{\bU_L\}$, consists of numerical solution trajectories computed by the numerical scheme on a low-resolution mesh. These low-fidelity trajectories are less accurate but more affordable compared to the high-fidelity ones. The key assumption behind the multi-fidelity architecture is that there exists a correlation between  data of different fidelity levels. In particular, given two consecutive solutions $\{\bU_L^m,\bU_L^{m+1} \}$ and $\{\bU_H^m,\bU_H^{m+1} \}$ with $\bU_L^m = \bU_H^m$, we assume that the relation between $\bU_L^{m+1}$ and $\bU_H^{m+1}$ is given by
\begin{equation}
    \label{eq:cor}
\bU^{m+1}_H = g(\bU^m_H, \bxi^{m},\bU^{m+1}_L),
\end{equation}
where the function $g$ takes the previous high-fidelity solution $\bU^m_H$, the normalized shifts $\bxi^m$, and the low-fidelity prediction $\bU_L^{m+1}$ as inputs, and
generates the high-fidelity solution $\bU_H^{m+1}$ with enhanced accuracy. 
While it is feasible to train a single-fidelity model using abundant high-fidelity data, where the model directly takes $\bU^m_H$ and  $\bxi^m$ as inputs and computes the approximation $\bU^{m+1}_H$ with high accuracy, in scenarios with limited availability of high-fidelity data, it is more effective to incorporate $\bU^{m+1}_L$ and explore the inherent correlation with \eqref{eq:cor}. Despite the low accuracy of the low-fidelity data $\bU^{m+1}_L$, it can still provide valuable information and greatly accelerate the training process of the high-fidelity component in the model.
By leveraging the inherent correlation between data of different fidelity levels, our multi-fidelity approach allows us to improve the accuracy of high-fidelity predictions even when high-fidelity data is limited. Moreover, by combining both high-fidelity and low-fidelity data in our training data set, we can benefit from the accuracy of the high-fidelity data and the accessibility of the low-fidelity data simultaneously. It enables us to develop a robust model that generalizes well across different fidelity levels and maximizes the utility of available data resources, leading to more efficient and accurate solutions in practice.

\begin{figure}[!htbp]
    \centerline{\includegraphics[width=0.95\textwidth]{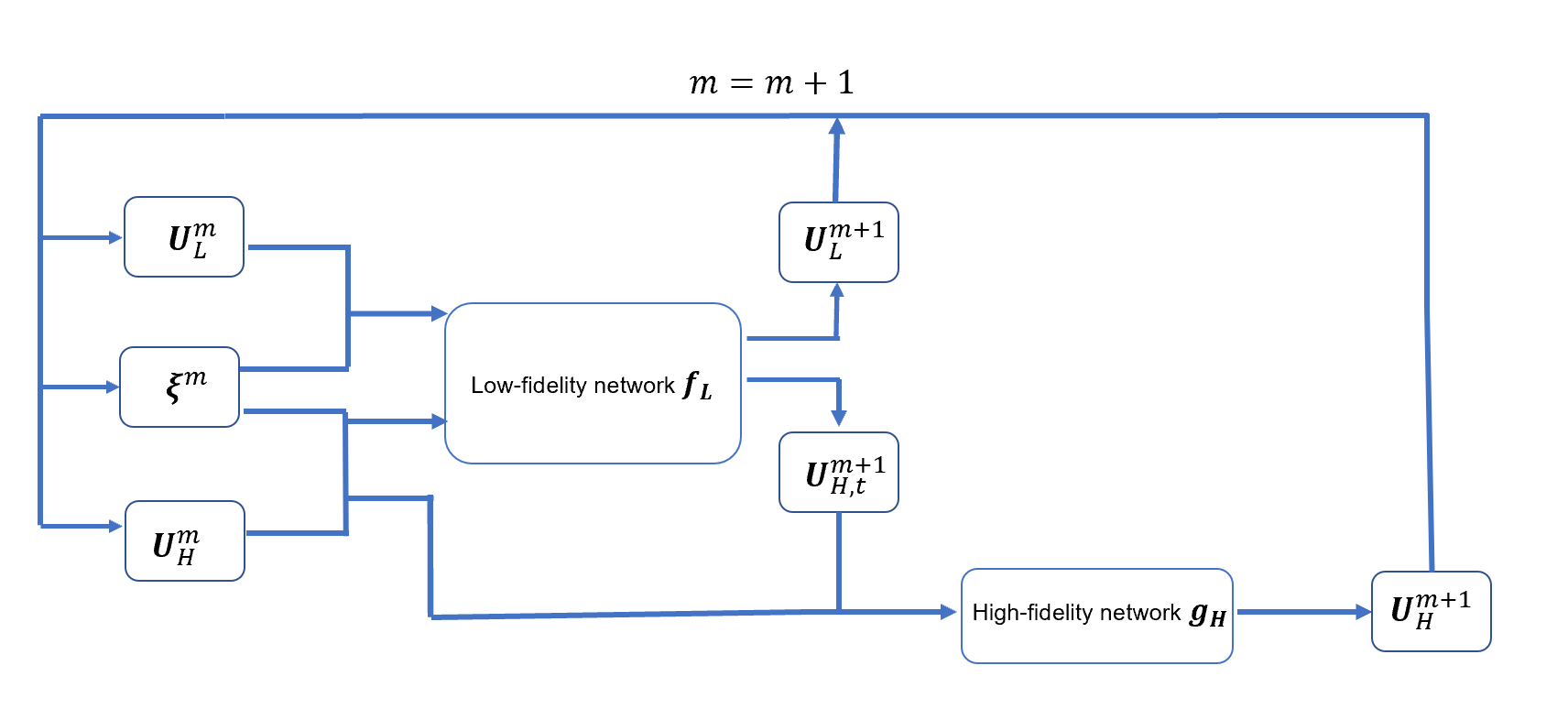}}
	\caption{Illustration of the multi-fidelity architecture.}\label{fig:multi_network_structure}
\end{figure}

Our multi-fidelity data-driven SL FV scheme employs a composite architecture, drawing inspiration from previous works under various frameworks such as Gaussian process regression \cite{raissi2016deep,perdikaris2017nonlinear}, PINNs \cite{meng_composite_2020}, 
and DeepONets \cite{howard_multifidelity_2022}.
Figure \ref{fig:multi_network_structure} provide a schematic diagram illustrating this architecture. It consists of two key components: the low-fidelity network, denoted as $f_L$, and the high-fidelity network, denoted as $g_H$. The low-fidelity network $f_L$ follows the same structure as the ML-based SL FV scheme described in Section \ref{sec:algorithm}, which takes a two-channel input ($\bU^m, \bxi^m$) and predicts the solution at the next time level, capturing the underlying trends in the data. Motivated by \eqref{eq:cor},  the high-fidelity network $g_H$ shares a similar structure to $f_L$, but it takes a three-channel input  ($\bU^m_H, \bxi^m, \bU^{m+1}_{H,t}$), where $\bU^{m+1}_{H,t}$ represents the intermediate solution computed by $f_L$ using the input ($\bU^m_H, \bxi^m$).   The objective of $g_H$ is to approximate the correlation and generate an enhanced solution $\bU^{m+1}_{H}$ with  an improved accuracy compared to $\bU^{m+1}_{H,t}$. By integrating both low-fidelity and high-fidelity components, our multi-fidelity architecture leverages the captured underlying trends from $f_L$ and incorporates the correlation approximation from $g_H$. This integration aims to enhance the accuracy of the solution beyond what each component can achieve individually. The composite architecture effectively combines the strengths of both components,  leading to improved accuracy and reliability in our data-driven SL FV scheme.

The training procedure for our multi-fidelity data-driven SL FV scheme involves the simultaneous training of the low-fidelity network $f_L$ and the high-fidelity network $g_H$. We have a set of $N_L$ low-fidelity trajectories $\left\{\{\bU^{(i),m}_L\}_{m=0}^{T_L}\right\}_{i=1}^{N_L}$ and $N_H$ high-fidelity trajectories $\left\{\{\bU^{(j),n}_{H}\}_{n=0}^{T_H}\right\}_{j=1}^{N_H}$ as training data, where
$\bU^{(i),m}_L$ represents the solution data in the $i$-th low-fidelity trajectory at time step $t_m$, and $\bU^{(j),n}_{H}$ represents the solution data in $j$-th high-fidelity  trajectory at  time step $t_n$. The initial conditions are drawn independently from the prescribed distribution.  The two networks, $f_L$ and $g_H$, are  trained by minimizing the one-step loss given by
\begin{equation}\label{eq:multi-loss}
    \mathcal{L}_{M F}\left(f_L,g_H\right)=\lambda_1 \mathcal{L}_{L F}+\lambda_2 \mathcal{L}_{H F},
\end{equation}
where  $\lambda_1$ and $\lambda_2$ are hyperparameters that balance low-fidelity loss $\mathcal{L}_{L F}$ and high-fidelity loss $\mathcal{L}_{H F}$. The low-fidelity loss  is defined as
\begin{equation}\label{eq:low_loss}
    \mathcal{L}_{L F}\left(f_L\right)=\frac{1}{N_L}\frac{1}{T_L} \sum_{i=1}^{N_L}\sum_{m=0}^{T_L-1}\left|f_L(\bU^{(i),m}_L, \bxi^{m})-\bU^{(i),m+1}_L\right|^2,
\end{equation}
where $f_L(\bU^{(i),m}_L, \bxi^{m})$ represents the prediction of the low-fidelity network $f_L$ for the next time step given the two-channel input ($\bU^{(i),m}_L$, $\bxi^{m}$).
The high-fidelity loss is defined as
\begin{equation}\label{eq:high_loss}
    \mathcal{L}_{H F}\left(g_H\right)=\frac{1}{N_H}\frac{1}{T_H} \sum_{j=1}^{N_H}\sum_{n=0}^{T_H-1}\left|g_H\left(\bU^{(j),n}_H,\bxi^{n},  \bU_{H,t}^{(j),n+1}\right)-\bU^{(j),n+1}_H\right|^2,
\end{equation}
where $\bU_{H,t}^{(j),n+1}=f_L\left(\bU^{(j),n}_H,\bxi^{n}\right)$ represents the intermediate solution computed by the the low-fidelity network $f_L$ using the two-channel input
$(\bU^{(j),n}_H,\bxi^{n})$, and $g_H\left(\bU^{(j),n}_H,\bxi^{n},  \bU_{H,t}^{(j),n+1}\right)$ represents the prediction of the high-fidelity network $g_H$ for the next time step given the three-channel input $\left(\bU^{(j),n}_H,\bxi^{n},  \bU_{H,t}^{(j),n+1}\right)$. During training, the objective is to minimize the multi-fidelity loss by adjusting the parameters of both networks. The hyperparameters $\lambda_1$ and $\lambda_2$ control the relative importance of the low-fidelity and high-fidelity losses in the overall loss function. Although it is possible to first train the network $f_L$ using the ML-based SL FV method discussed in Section \ref{sec:algorithm} and subsequently train the network $g_H$ with fixed $f_L$ by \eqref{eq:multi-loss}, simultaneous training of the neural networks with multi-task learning offers several benefits, including improved data efficiency and accelerated learning, as highlighted in \cite{crawshaw2020multi}. 
It is worth noting that one-step training, as described above, may suffer slightly weaker generalization. On the other hand, unrolling the training over multiple time steps can improve the accuracy and stability at the cost of increased training difficulty, as discussed in \cite{brandstetter2022message}. 
In the conducted numerical experiments, one-step training is used for linear cases, while training is unrolled for 16 time steps when dealing with the nonlinear VP system.

For the implementation of the trained multi-fidelity SL FV scheme in the testing phase,
the solutions are updated in time as follows: Starting from the solution $\bU^m$ at time level $m$, where $\bU^0$ is initialized as the cell averages of the given initial condition, the characteristic equation \eqref{eq:characteristic} is solved to obtain the normalized shifts $\bxi^m$. Next, using  the trained low-fidelity network $f_L$, the intermediate solution $\bU^{m+1}_{t}$ is computed by taking the two-channel inputs $(\bU^m,\bxi^m)$, i.e., $\bU^{m+1}_{t} = f_L(\bU^m,\bxi^m)$.   Finally, 
the solution $\bU^{m+1}$ with enhanced accuracy is generated by feeding $\bU^m$, $\bxi^m$, and $\bU^{m+1}_{t}$ as inputs to the trained high-fidelity network $g_H$, resulting in $\bU^{m+1} = g_H(\bU^m,\bxi^m,\bU^{m+1}_{t})$.
It is numerically observed that the intermediate solution $\bU_{t}^{m+1}$ achieves a higher level of accuracy compared to $\bU_{L}^{m+1}$, which is obtained by the single-fidelity 
SL FV method trained solely with low-fidelity data. This difference in accuracy validates the effectiveness of the multi-fidelity approach and the accelerated learning process for the high-fidelity network $g_H$, which benefits from the more accurate intermediate solution as part of its inputs.


\subsection{The Vlasov-Poisson system}\label{sec:vp}
In this subsection, we extend the proposed algorithm in Section \ref{sec:scheme:transport} for linear transport equations to the nonlinear VP system, to address the challenges posed by  its inherent nonlinearity.

The VP system is a fundamental model in plasma physics that describes interactions between charged particles through self-consistent electric fields, modeled by Poisson’s equations in the non-relativistic zero-magnetic field limit. The dimensionless governing equations under 1D in space and 1D in velocity (1D1V) setting is given by
\begin{align}
    &f_t + vf_x + E(x,t)f_v = 0,\quad (x,v)\in\Omega_x\times\Omega_v,\label{eq:vlasov}\\
    &E_x = \rho - 1,\quad \rho(x,t) = \int_{\Omega_v} f(x,v,t) dv,\label{eq:poisson}
\end{align}
where $f(x,v,t)$ is the probability distribution function of electrons at position $x$ with velocity $v$ at time $t$, and  $\rho$ denotes the density. $\Omega_x$ denotes the physical domain, while  $\Omega_v$ represents the velocity domain. Here we assume a uniform background of fixed ions under a self-consistent electrostatic field $E$.

Unlike the linear transport equation \eqref{eq:transport1d}, the Vlasov equation \eqref{eq:vlasov} is a nonlinear transport equation, which introduces additional challenges for accurately approximating its characteristic equations
\begin{equation}\label{eq:vpchar}
\begin{cases}
 \displaystyle  \frac{d x(t)}{d t} &= v(t),\\[3mm]
  \displaystyle  \frac{d v(t)}{d t} &= E(x(t),t).
\end{cases}
\end{equation}
Consequently, designing an SL scheme becomes more complex for the VP system.
To address these challenges, the splitting approach has gained popularity in the literature \cite{cheng1976integration,sonnendrucker1999semi}
. This approach decouples the VP system into several linear transport equations, which are much easier to solve in the SL framework. However, it is important to note that the inherent splitting error  can potentially compromise the overall accuracy of the solution. 
Recently, a non-splitting SL methodology has been introduced for the VP system, utilizing  the commutator-free RKEI \cite{celledoni2003commutator}. 
By employing RKEI, the VP system is decomposed into a sequence of linearized transport equations with frozen coefficients \cite{cai_high_2021,zheng2022fourth}, 
which can be effectively solved using the proposed multi-fidelity ML-based SL FV scheme. Our novel approach combines the advantages of RKEI and multi-fidelity ML-based SL FV schemes, thereby avoiding the splitting errors inherent in traditional splitting approaches and enabling efficient and accurate solutions to the VP system.

To introduce the algorithm, we begin by considering a uniform partition of the domain $\Omega_x\times\Omega_v$ with $N_x\times N_v$ cells, i.e., $\Omega_x\times\Omega_v=\bigcup_{ij}I_{ij}$. The mesh sizes in the $x$ and $v$ directions are represented by as $h_x$ and $h_v$, respectively. Let  $f^m_{ij}$ denote the numerical approximation of the cell average  of $f$ in the cell $I_{ij}$ at time level $t^m$ and $E_{i+\frac{1}{2}}^m$ be the numerical solution of $E$ at the endpoint of the cell in the $x$ direction. The collections $\{E^m_{i-\frac12}\}$ and $\{f^m_{ij}\}$ are denoted by $\bE^m$ and $\bF^m$, respectively.
Within the FV framework,  
the electric field $E(x)$ can be straightforwardly approximated from Poisson's equation \eqref{eq:poisson} as
\begin{equation}
\label{eq:elec}
E^m_{i+\frac12} = E^m_{-\frac12} + h_x\sum_{l=1}^i\left(h_v\sum_{j=1}^{N_v}f_{ij}^m-1\right),\quad i=1,\dots,N_x
\end{equation}
Here, $E^m_{-\frac12}$ is determined by the given boundary conditions. For instance, in case of periodic boundary conditions, we have $E^m_{-\frac12}=0$.

\begin{table}[!htbp]
\begin{minipage}{.5\linewidth}
    \centering

    \caption{First-order RKEI}
    \label{tab:1strkei}

    \medskip

\begin{tabular}{l|l}
0 & 0 \\
\hline & 1
\end{tabular}
\end{minipage}\hfill
\begin{minipage}{.5\linewidth}
    \centering

    \caption{Second-order RKEI}
    \label{tab:2ndrkei}

    \medskip

\begin{tabular}{l|ll}
0 & & \\
$\frac{1}{2}$ & $\frac{1}{2}$ & 0 \\
\hline & 0 & 1
\end{tabular}
\end{minipage}
\end{table}

Our approach employs the RKEI technique to alleviate the difficulties associated with accurately tracking the characteristics for the nonlinear VP system. The RKEI method is represented by a Butcher tableau, which specifies the coefficients for the integration steps. The accuracy of the method is determined by order conditions. The simplest RKEI method is given by the Butcher tableau \ref{tab:1strkei}, which is first order accurate \cite{celledoni2003commutator}.
With this first-order RKEI, we can develop a  ML-based SL FV scheme for the VP system as follows: 
\begin{itemize}
    \item[(1)] Compute $\bE^m$ using  \eqref{eq:elec} based on the cell average approximations $\bF^m$. 
    \item[(2)] Linearize the Vlasov equation with the fixed electric field  $\bE^m$ and solve the characteristic equations \eqref{eq:vpchar}  to obtain the normalized shifts.
    \item[(3)] Use the ML-based SL FV method, together with $\bF^m$ and the normalized shifts, to predict $\bF^{m+1}$.
\end{itemize}
The procedure  can be summarized as
\begin{equation}
\label{eq:1stvp}
\bF^{m+1} = MLSL(\bE^m,\Delta t)\bF^{m},
\end{equation}
where $MLSL(\bE^m,\Delta t)$ denotes the ML-based SL FV evolution operator for the Vlasov equation with the fixed electric field  $\bE^m$ and time step $\Delta t$.
However, the low order temporal accuracy of the first-order RKEI  can limit its  generalization capability.  

In the following, we enhance the performance by employing a second order RKEI \cite{celledoni2003commutator},
represented by the Butcher tableau \ref{tab:2ndrkei}.
The corresponding ML-based SL FV algorithm for the VP system can be summarized as 
\begin{equation}
\label{eq:2ndvp}
\begin{cases}
        \bF^{m,*} = MLSL(\bE^m,\frac12\Delta t)\bF^{m},  \\
    \bF^{m+1} = MLSL(\bE^{m,*},\Delta t)\bF^{m}, 
\end{cases} 
\end{equation}
where $\bE^m$ and $\bE^{m,*}$ are determined via  \eqref{eq:elec} by $\bF^{m}$ and $\bF^{m,*}$, respectively. The second order scheme involves an intermediate stage $\bF^{m,*}$ and requires two applications of the ML-based SL FV evolution operator. It is numerically demonstrated that this second order method \eqref{eq:2ndvp} exhibits improved  generalization capabilities and achieves a higher level of accuracy compared to the first counterpart \eqref{eq:1stvp}. By combining the proposed multi-fidelity framework with the ML-based SL approach, we can effectively learn an optimal discretization for the VP system, while significantly reducing the dependence on high-fidelity data. Moreover, all the desired properties of the method for solving linear equations, such as mass conservation and translational equivariance, are preserved for the nonlinear VP system. The learned multi-fidelity model  achieves a level of accuracy that surpasses traditional numerical algorithms with comparable mesh resolution, leading to significant computational savings.

\section{Numerical results}
\label{sec:num}
In this section, we carry out a series of numerical experiments to demonstrate the performance of the proposed multi-fidelity ML-assisted SL FV scheme for various benchmark 1D and 2D linear transport equations, as well as 1D1V nonlinear VP systems.  Note that the performance of the proposed scheme is significantly influenced by the choice of hyperparameters for the NN structure. Here, we present numerical results using the following default settings. For the linear equations, we utilize 6 convolutional layers with 32 filters each, and a kernel size of 5 for both 1D and 2D cases. For the nonlinear VP systems, we employ 9 convolutional layers with 32 filters each, using a kernel size of 5. In addition,  to ensure mass conservation, a constraint layer is added.  Throughout the training process, we employ the ReLU function  as the activation function and use the Adam optimization algorithm \cite{Ilya_fix_2017}.

We are allowed to employ any accurate and stable FV schemes to generate both high-fidelity and low-fidelity data. Specifically, for linear transport equations, we adopt the Eulerian fifth-order FV WENO (WENO5) method 
\cite{jiang1996efficient}, 
coupled with the third-order strong-stability-preserving RK  time integrator \cite{gottlieb2001strong}.
For nonlinear VP systems, we employ the conservative SL FV WENO5 scheme coupled with a fourth-order RKEI 
\cite{zheng2022fourth}.  
The high-fidelity training data is obtained by coarsening high-resolution solution trajectories into a low resolution mesh. The low-fidelity training data is generated by directly running the transport scheme on the low-resolution mesh. In all test examples, we primarily report the results of the proposed multi-fidelity learning approach and compare these results with the single-fidelity model developed in Section \ref{sec:algorithm} and the traditional WENO5 scheme. 
The plots presented below indicated the performance of each approach: ``Multi-fidelity" refers to the proposed method, ``High-fidelity" denotes the single-fidelity learned SL FV  method trained solely with high-fidelity data, and ``Low-fidelity" denotes the learned SL FV  method trained solely with low-fidelity data. 

\subsection{Transport equations}
In this subsection, we present numerical results for simulating several 1D transport equations.

\begin{example}\label{square}
In this example, we consider the following advection equation with a constant coefficient
\begin{equation}\label{eq:const1d}
u_t + u_x = 0,\quad x\in [0,1],
\end{equation}
and periodic conditions are imposed. 
\end{example}

We generate 30 low-resolution solution trajectories with a 32-cell grid as low-fidelity training data. Each trajectory comprises 54 sequential time steps. For high-fidelity training data, we downsample 15 high-resolution solution trajectories computed with a 256-cell grid by a factor of 8. Each trajectory consists of 15 sequential time steps. Moreover, all trajectories, both low- and high-fidelity, are initialized a square wave with the height randomly sampled from the interval $[0.1, 1]$ and the width randomly selected from $[0.2, 0.4]$.
 We use a CFL number of 1.8 for determining the time step and a centered 5-cell stencil for updating the solution. During the multi-fidelity model training, the parameters $\lambda_1$ and $\lambda_2$ are set as $0.1$ and $1$, respectively. To test, we randomly generate square functions with  widths and heights within the same range as initial conditions.
 For comparison, we generate a reference solution following the same approach as  the high-fidelity data creation. 
 
\begin{figure}[!htbp]
 \centering
    \centerline{\includegraphics[width=0.95\textwidth]{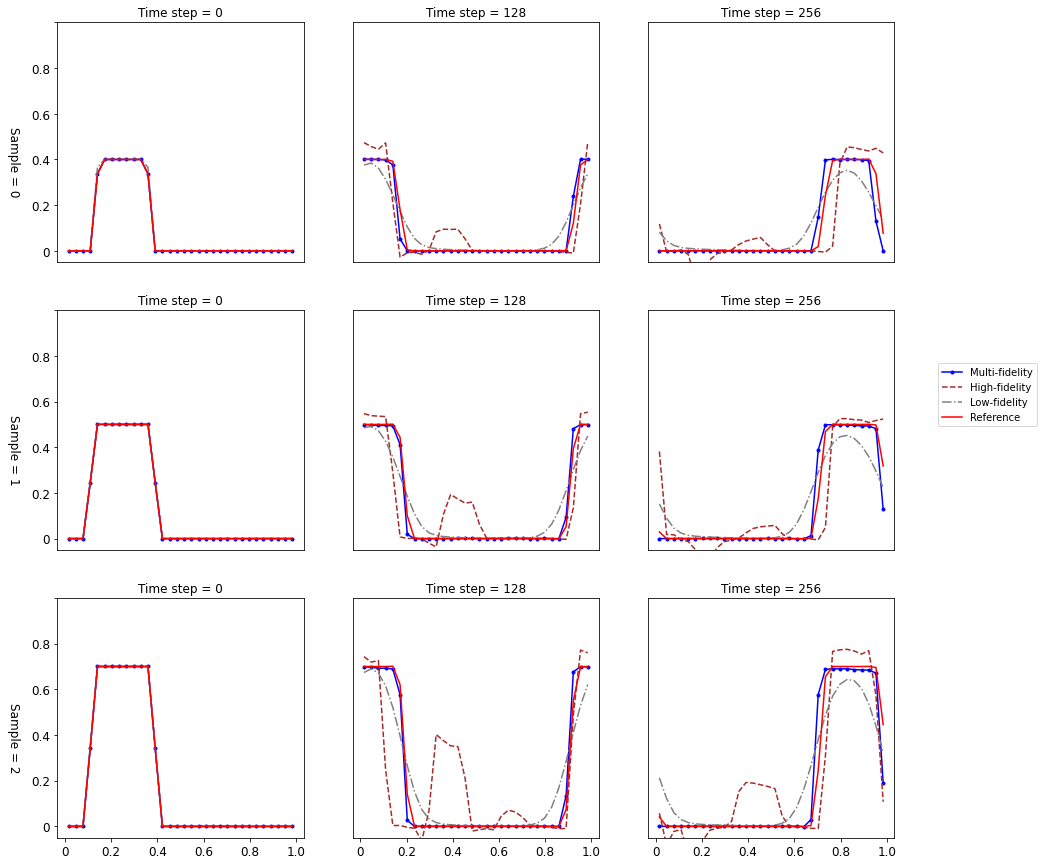}}
	\caption{Three test samples for square waves in Example \ref{square}.  }\label{fig:square_three_curve}
\end{figure}

Figure \ref{fig:square_three_curve} displays three test samples obtained from forward integration at different time instances, computed by the proposed multi-fidelity model, as well as by single-fidelity models trained solely with high-fidelity or low-fidelity data. It is observed that the low-fidelity model exhibits significant smearing near discontinuities due to the inherent low accuracy of the training data. Meanwhile, the high-fidelity model generates highly inaccurate results due to the lack of sufficient high-ﬁdelity data. In contrast, our proposed multi-fidelity approach demonstrates superior shock resolution, with sharp shock transitions and no spurious oscillations. Notably, even after a long-time simulation of 512 time steps, which is about 34 times the number of time steps used for generating high-fidelity training data, our multi-fidelity method remains remarkably close to the reference solution. This validates the high stability and accuracy of  the learned discretization. 

\begin{figure}[!htbp]
    \centerline{\includegraphics[width=0.95\textwidth]{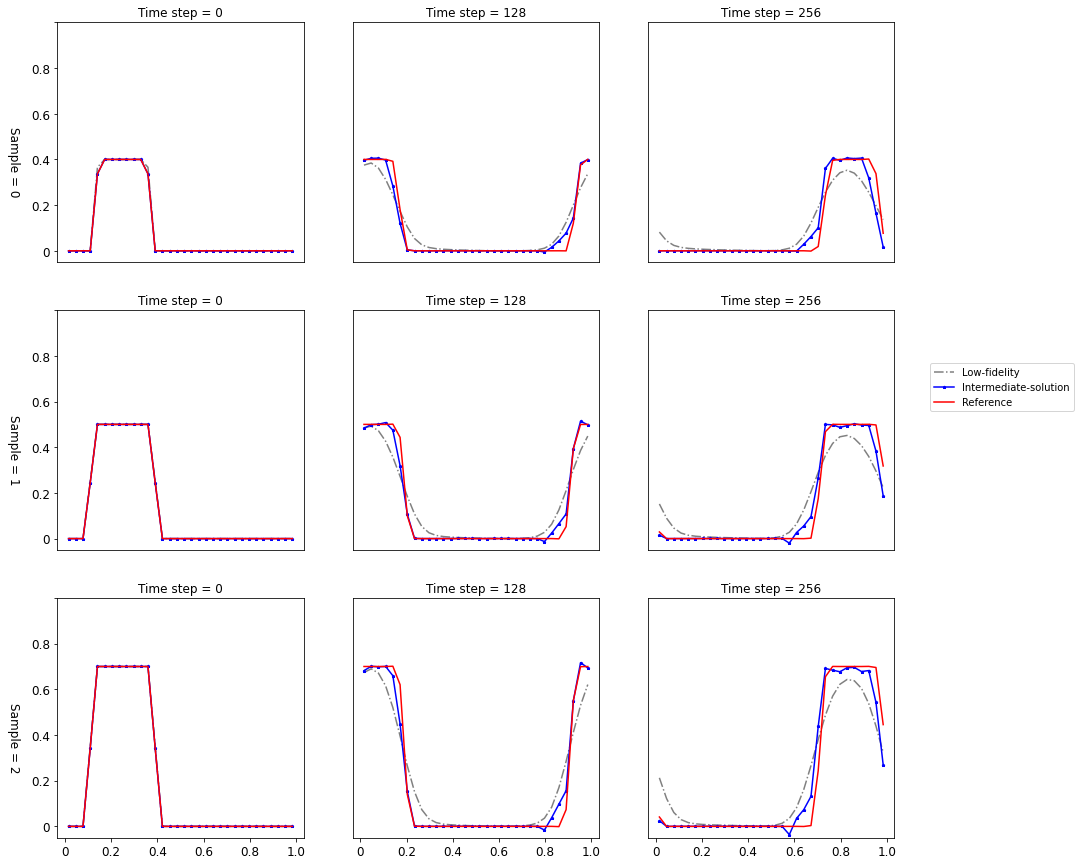}}
	\caption{The intermediate solutions of three test samples for square waves in Example \ref{square}.  }\label{fig:square_inter}
\end{figure}

\begin{figure}[!htbp]    
 \centering
 \subfigure[Example \ref{square}]{\includegraphics[width=0.4\textwidth]{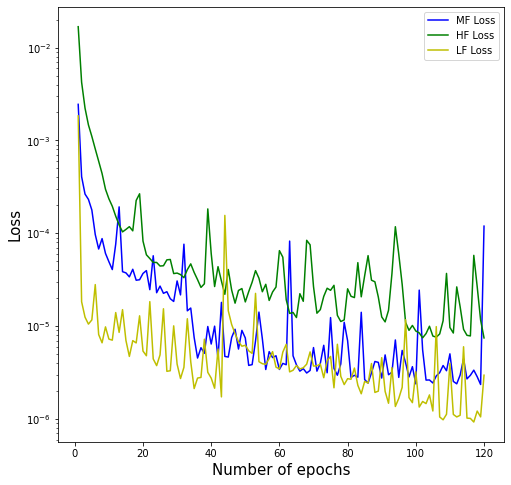}}
 \quad\quad\quad
 \subfigure[Example \ref{triangle}]{\includegraphics[width=0.4\textwidth]{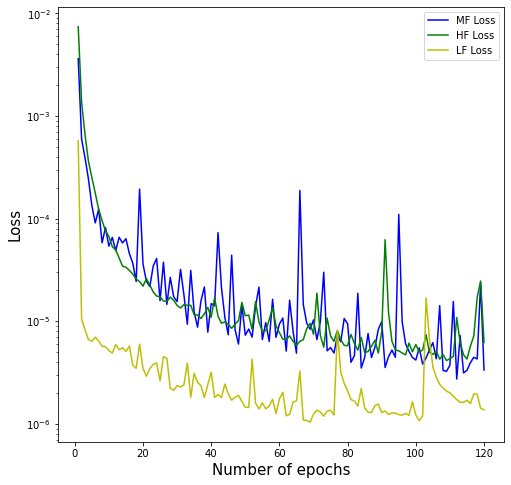}}
	\caption{The training loss for three different methods in Examples \ref{square} (left) and  \ref{triangle} (right). }
 \label{fig:loss}
\end{figure}

Figure \ref{fig:square_inter} shows the intermediate solution of our multi-fidelity model for the three test samples, obtained from the trained low-fidelity network $f_L$ of the proposed multi-fidelity method. It is observed that the intermediate solution is significantly more accurate compared to the results obtained by the low-fidelity model alone. This indicates that the simultaneous training of the multi-fidelity model is beneficial for accelerating the learning process of the high-fidelity component. For brevity, we omit the comparison result between the proposed ML-based method and the traditional  WENO5 scheme. However, it is worth mentioning that the results obtained with WENO5 are similar to those of the low-fidelity model. Figure \ref{fig:loss}(a) illustrates the training loss of the three models throughout their respective training epochs.

\begin{example}\label{triangle}
In this example, we consider the advection equation \eqref{eq:const1d} with a more complicated solution profile consisting of triangle and square waves.
\end{example}

Similar to the previous example, we generate 30 low-resolution solution trajectories using a 32-cell grid  as the low-fidelity training data. Each trajectory consists of 54 sequential time steps. To create high-fidelity training data, we downsample 15 trajectories of high-resolution solutions computed with a 256-cell grid by a factor of 8. Each trajectory contains 15 sequential time steps. For both low- and high-fidelity trajectories, the initial conditions include one triangle centered at 0.25, and one square wave centered at 0.75. The heights and widths of these shapes are randomly sampled from  the ranges [0.2, 0.8] and  [0.2, 0.3], respectively. We use a fixed stencil size of  3 and set the CFL number to 1.8. The parameters for training are chosen as $\lambda_1=0.05$ and $\lambda_2=1$. For testing, the initial conditions are randomly generated from the same distribution used to create the training data. 

\begin{figure}[!htbp]
    \centerline{\includegraphics[width=0.95\textwidth]{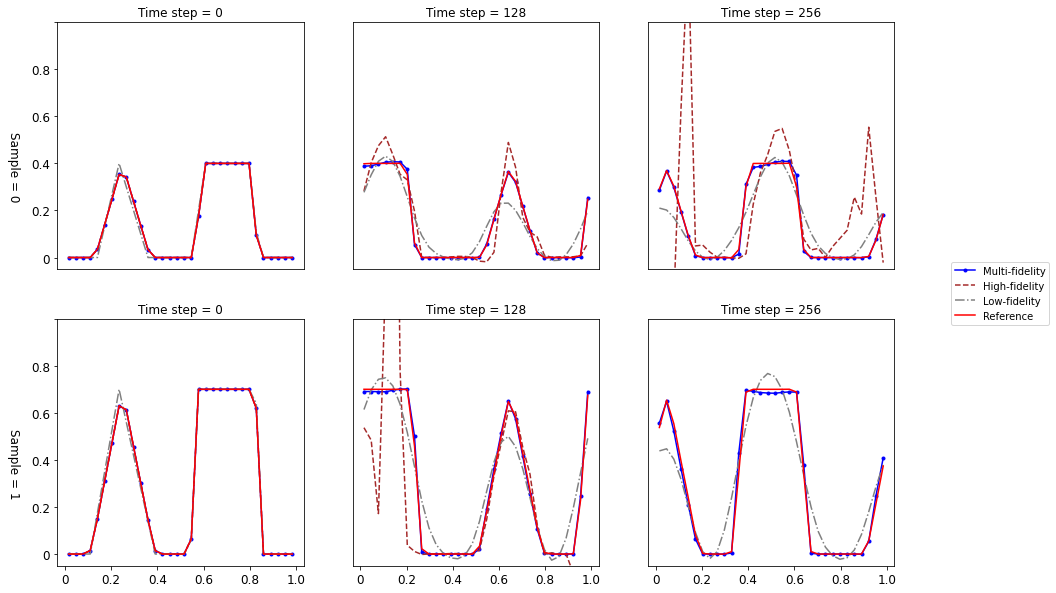}}
	\caption{Numerical solutions of two test samples for advection of triangle and square waves in Example \ref{triangle}.  CFL=1.8.}\label{fig:triangle_three_curve}
\end{figure}

\begin{figure}[!htbp]
    \centerline{\includegraphics[width=0.95\textwidth]{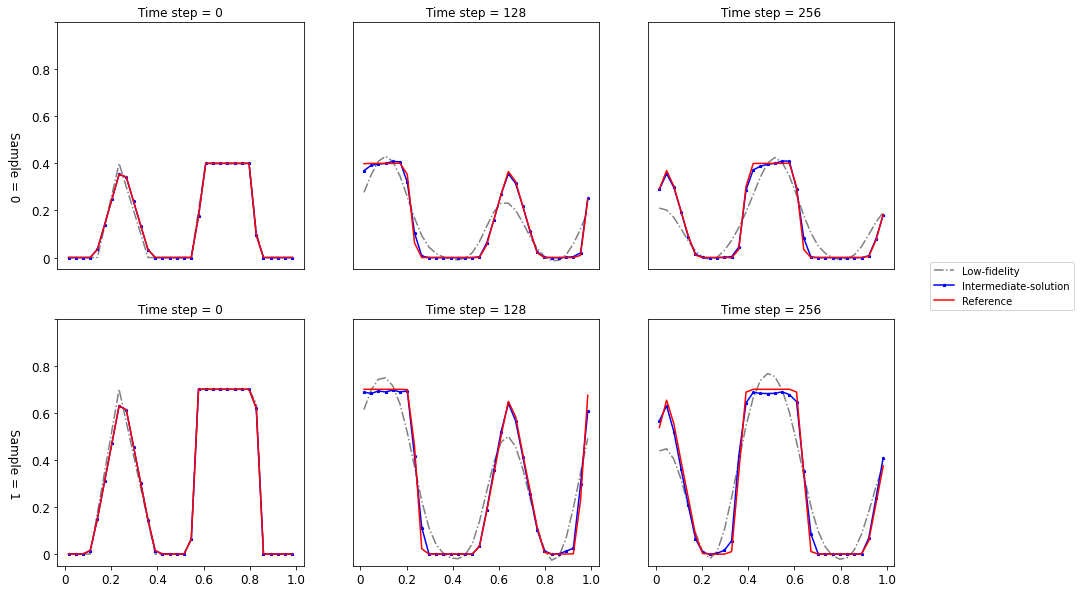}}
	\caption{The intermediate solutions of two test samples for advection of triangle and square waves in Example \ref{triangle}.  CFL=1.8.}\label{fig:triangle_inter}
\end{figure}

In Figure \ref{fig:triangle_three_curve}, we plot two test samples at several time instances  during forward integration. It is observed that the proposed multi-fidelity solver generates numerical results with  significantly higher resolution of the underlying non-smooth structures compared to the low-fidelity model. Moreover, the high-fidelity model exhibits severe spurious oscillations and eventually blows up. Figure \ref{fig:triangle_inter} presents the intermediate solution of our multi-fidelity method for the two test samples. Similar to the previous example, we observe that the intermediate solution of our multi-fidelity method is much more accurate compared to the results obtained by the low-fidelity model. The training loss for three different methods during their respective training is shown in Figure \ref{fig:loss}(b).
 


\begin{example}
\label{eg:variable}
In this example, we simulate the following 1D advection equation with a variable coefficient
\begin{equation}\label{eq:variable-velocity}
    u_t+(u\sin(x+t))_x=0,\quad x\in[0,2\pi],
\end{equation}
subject to periodic conditions. This example is more challenging than the previous two examples. In addition to the pure shift, the solution profiles will gradually deform over time and exhibit more complicated structures.
\end{example} 

To generate low-fidelity training data, we produce 90 solution trajectories on a 32-cell grid, with each trajectory consisting of 4 sequential time steps. For high-fidelity training data, we coarsen 90 solution trajectories over a 256-cell grid by a factor of 8, with each trajectory consisting of 4 sequential time steps.  As in the first example, the initial condition is a step function with heights randomly sampled from the range $[0.1,1]$
 and widths sampled from the range $[2.5,3.5]$. The center of each step function is randomly selected from the entire domain of 
$[0,2\pi]$. The stencil size is set to be 5, and the CFL number is 0.5. 
 We choose the loss weights as $\lambda_1=1$ and $\lambda_2=1$.

\begin{figure}[!htbp]
    \centerline{\includegraphics[width=0.9\textwidth]{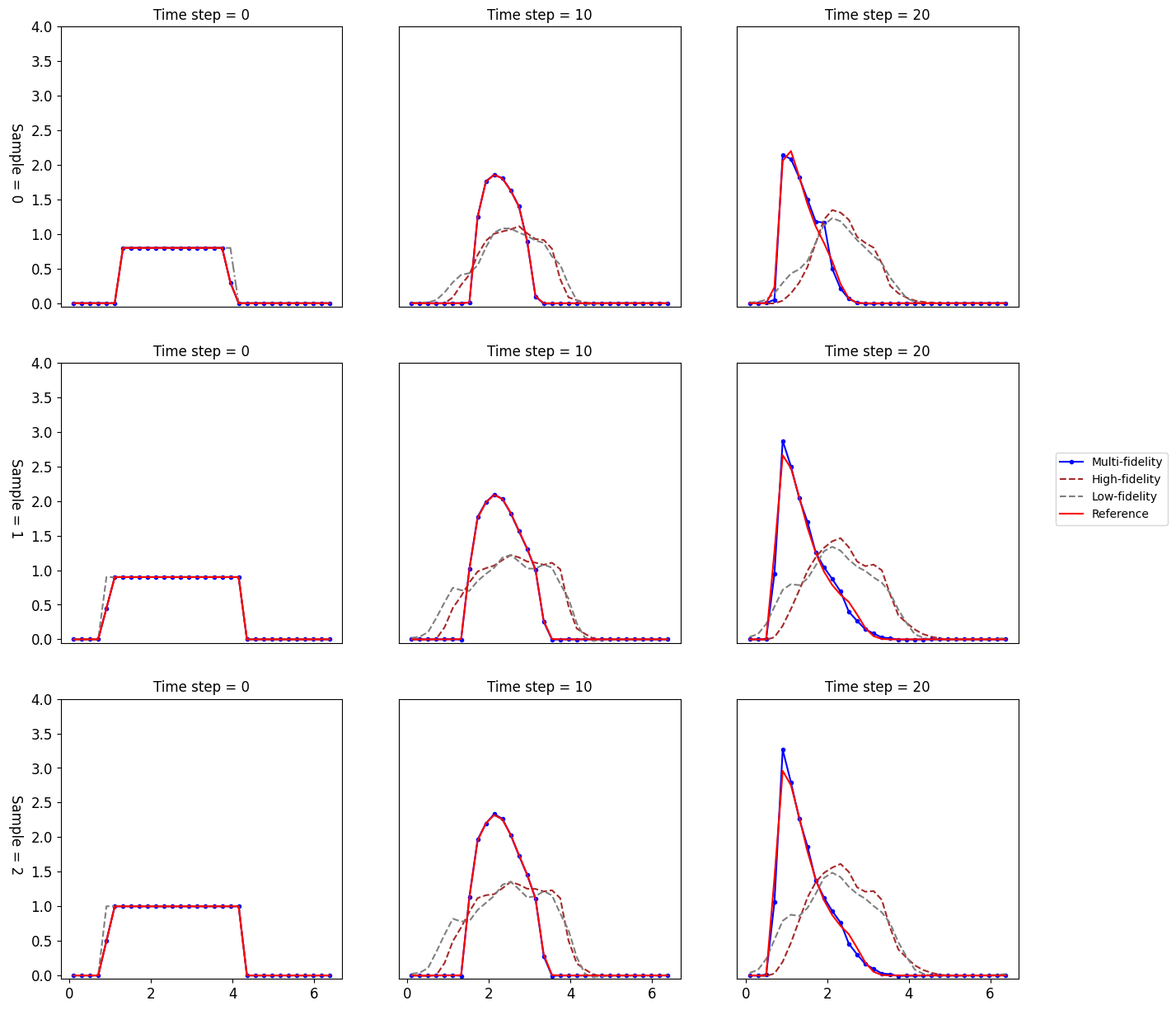}}
	\caption{Numerical solutions of three test samples for the transport equation with a variable coefficient in Example \ref{eg:variable}. CFL=0.5.}\label{fig:variable_three_curve}
\end{figure}

Figure \ref{fig:variable_three_curve} shows three test samples plotted at various time instances during the forward integration. We run the trained models up to 20 time steps, which is greatly beyond the range of the training data, to assess their generalization capability. Our multi-fidelity model is able to accurately resolve such highly deformed solution structures that are not seen in the training data set. In contrast, both the high-fidelity and low-fidelity solvers fail generalize well and produce irrelevant results.


\begin{example}\label{defor}
In this example, we simulate the 2D deformational flow problem proposed in 
\cite{leveque_high-resolution_1996},
governed by the following  transport equation
  \begin{equation}\label{eq:2d:defor}
u_t + (a(x,y,t)u)_x+(b(x,y,t)u)_y = 0,\quad (x,y)\in[0,1]^2, 
\end{equation}
with the velocity field being a periodic swirling flow \begin{equation}\label{eq:deformational flow}
\begin{aligned}
    a(x,y,t) &= \sin^2(\pi x)\sin(2\pi y)\cos(\pi t/T),\\
    b(x,y,t) &= -\sin^2(\pi y)\sin(2\pi x)\cos(\pi t/T),
\end{aligned}
\end{equation}
where $T$ is a constant.
This is a widely recognized benchmark test for numerical transport solvers. It exhibits a distinct dynamics, in which the solution profile deforms over time as it follows the flow. The direction of the flow reverses at $t=T/2$, and the solution returns to its initial state at $t=T$, completing a full cycle of the evolution.
\end{example}
We set  $T=2$ and let the initial condition be a cosine bell centered at $[r_x,r_y]$:
\begin{equation}\label{eq:deformation-initial}
\begin{aligned}
    u(x,y) &= \frac{1}{2}[1+\cos (\pi r)]\\
    r(x,y) &= \min\left[1,6\sqrt{(x-r_x)^2+(y-r_y)^2}\right].
\end{aligned} 
\end{equation} 
We generate 4 trajectories using a high-resolution mesh of $256\times256$ cells, which are subsequently coarsened by a factor of 8 in each dimension to create the high-fidelity training data on the mesh of $32\times 32$ cells. Each trajectory  contains a sequence of 214 time steps from $t=0$ to $t=T$. To obtain the low-fidelity training data, we generate 18 trajectories on the  mesh of $32\times 32$ cells. Each trajectory also contains a sequence of 214 time steps from $t=0$ to $t=T$. The initial conditions for all trajectories are given by \eqref{eq:deformation-initial}, where $r_x$ and $r_y$ are randomly sampled from the range $[0.25,0.75]$. We set the CFL number to 1.8 and use a stencil of size $5\times5$.   For training setting, the loss weights are chosen as $\lambda_1=0.1$ and $\lambda_2=1$.

\begin{figure}[!htbp]
    \centerline{\includegraphics[width=0.9\textwidth]{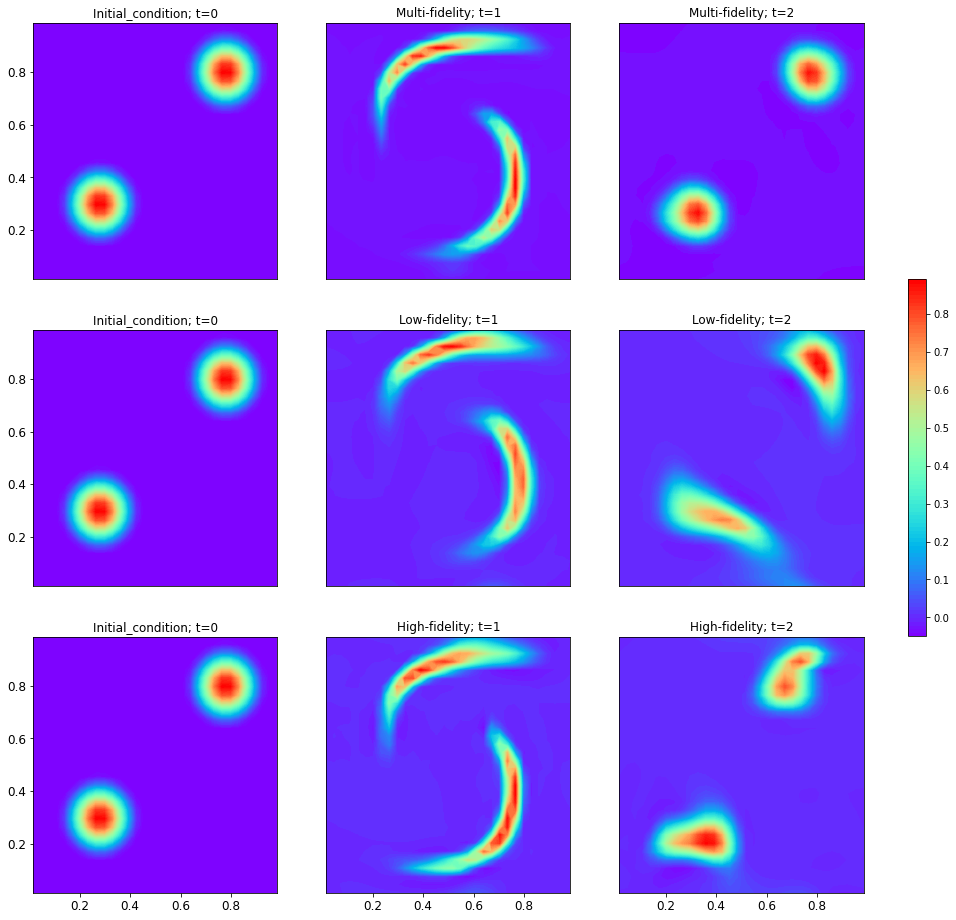}}
	\caption{Contour plots of the numerical solutions for the 2D deformational flow at $t=0,1,2$ in Example \ref{defor}. }\label{fig:two_deform}
\end{figure}

Although the training data consists of solution trajectories featuring a single bell, we demonstrate that the trained multi-fidelity model can generalize to simulate problems with an initial condition that contains two randomly placed bells.
In Figure \ref{fig:two_deform}, we show contour plots of numerical solutions computed by the proposed multi-fidelity method, the high-fidelity model, and the low-fidelity model. 
It is observed that the proposed method greatly surpasses single-fidelity models in capturing deformations and accurately recovering the initial profile at $t=T$.
 
 Furthermore, we demonstrate the efficiency of the proposed method by providing the comparison of the error and run-time between the proposed multi-fidelity model and the traditional Eulerian RK WENO5 method.  In Table
 \ref{table:defor-time} we provide the run-time for simulating three test samples up to $t=T$ using the multi-fidelity model with a mesh resoluton of $32\times32$ cells and CFL number of 1.8, as well as the RK WENO5 with mesh resolutions of  $32\times32$ cells  and   $128\times128$ cells, both using the CFL number of 0.6. Table \ref{table:defor-error} reports the corresponding mean square errors. 
The computational time of the proposed multi-fidelity model is slightly higher than that of the RK WENO5 method with the same mesh resolution of $32\times 32$ cells. However, it is significantly lower than the RK WENO5 over a finer mesh of $128\times 128$ cells. Meanwhile, the multi-fidelity model achieves much smaller errors compared to the RK WENO5 method with the same mesh size, and only slightly larger errors than the RK WENO5 method using the high-resolution mesh.

\begin{table}[!htbp]
		\centering
		\caption{Run-time comparison of Example \ref{defor}. Run-time (seconds) measure for one period from $t=0$ to $t=T$ on a GeForceRTX 4070 Ti GPU for the proposed multi-fidelity model, and on a CPU for the traditional Eulerian RK WENO5 method using Python.} \label{table:defor-time}
		\smallskip
		\begin{tabular}[c]{c c c c  }
			\hline
			 Samples\textbackslash Method & Multi-fidelity($32\times 32$) &  WENO5($32\times 32$) & WENO5($128\times 128$) \\
			\hline
		Sample = 0	    &  3.5734 & 1.3498 & 26.9169   \\  \hline
            Sample = 1 &  3.4712 &1.3766 &27.0657 \\
			\hline
            Sample = 2 &  3.3923 &1.3200 &26.1897 \\
			\hline

		\end{tabular}
	\end{table}

\begin{table}[!htbp]
		\centering
		\caption{Mean square errors of Example \ref{defor}.} \label{table:defor-error}
		\smallskip
		\begin{tabular}[c]{c c c c  }
			\hline
			Samples \textbackslash Method & Multi-fidelity($32\times 32$) &  WENO5($32\times 32$) & WENO5($128\times 128$) \\
			\hline
           Sample = 0 & 1.0249E-5 &2.0531E-3 & 5.8772E-6 \\
			\hline
             Sample = 1 & 1.8541E-5 &1.6618E-3 & 5.1721E-6 \\
			\hline
            Sample = 2 & 1.6547E-5 &1.6547E-3 & 5.2394E-6 \\
			\hline

		\end{tabular}
	\end{table}

\subsection{Nonlinear Vlasov-Possion System}
In this subsection, we present the numerical results for simulating the nonlinear 1D1V VP system. We demonstrate the efficiency and accuracy of the proposed multi-fidelity SL FV method coupled with the second-order RKEI by comparing it to single-fidelity models and the traditional SL FV WENO5 scheme. Additionally, we numerically verify the advantage of the second-order RKEI over the first-order scheme. The training data and reference solutions are generated by the fourth-order conservative SL FV WENO scheme 
\cite{zheng2022fourth}.

\begin{example}\label{landau-damping}
In this example, we consider the Landau damping with the initial condition 
    \begin{equation}\label{eq:landau-damping}
        f(x,v,t=0)=\frac{1}{\sqrt{2\pi}}(1+\alpha \cos(kx))\exp\left(-\frac{v^2}{2}\right), \quad x\in[0,L], \quad v\in[-V_c,V_c],
    \end{equation}
    where $k=0.5$, $L=4\pi$, and $V_c=2\pi$. 
\end{example}

We generate 6 solution trajectories with a $32\times 64$ grid to obtain the low-fidelity training data. The high-fidelity training data is generated by coarsening 2 high-resolution solution trajectories on a $256\times 512$-cell grid by a factor of 8 in each dimension. The initial conditions for both the low- and high-fidelity training data sets are determined using \eqref{eq:landau-damping}, with  $\alpha$ randomly selected from a uniform distribution in the range $[0.05,0.45]$. We set the stencil size to be $5\times 5$ and maintain a chosen CFL number of 1.8. For training setting, the loss weights are chosen as $\lambda_1=0.1$ and $\lambda_2=1$.
For the purpose of comparison, the reference solution is generated with the same approach used to create the high-fidelity data. 

\begin{figure}[!htbp]
\centering
    \centerline{\includegraphics[width=0.9\textwidth]{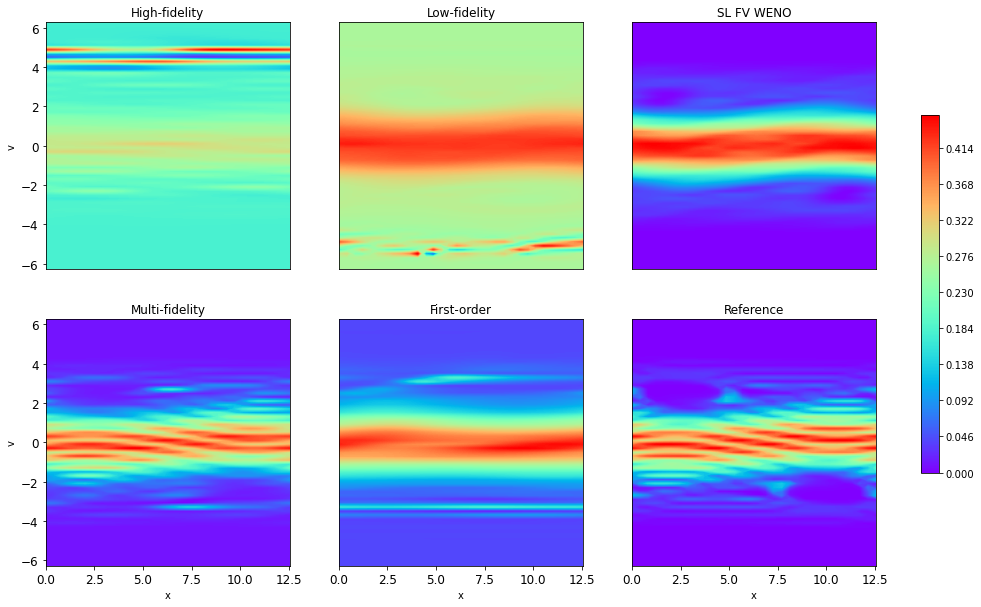}}
	\caption{Contour plots of numerical solutions of the strong Landau damping at $t=40$ in Example \ref{landau-damping} with $\alpha=0.5$.  ``Multi-fidelity" denotes our multi-fidelity method coupled with the second-order RKEI. ``First-order" denotes our multi-fidelity method with the first-order RKEI.}\label{fig:landau_contour}
\end{figure}

For testing, we set $\alpha=0.5$, yielding the strong Landau damping, which lies outside the range of the training data. Figure \ref{fig:landau_contour} presents contour plots of the numerical solutions computed by our multi-fidelity method with the second order RKEI, high-fidelity solver, low-fidelity solver, the multi-fidelity scheme with the first-order RKEI, and the traditional SL FV WENO method, all implemented on the same mesh with  $32\times 64$ cells.  It can be observed that the proposed method with the second-order RKEI captures the filamentation structure, while both the single-fidelity models generate inaccurate results. In addition, the multi-fidelity scheme coupled with the first-order RKEI fails to capture the fine-scale structures of interest. This limitation arises due to the low order accuracy in time, which restricts its ability to generalize effectively. The traditional SL FV WENO scheme produces reasonable results but exhibits smeared solution structures, primarily due to the low mesh resolution.

\begin{figure}[!htbp]
    \centerline{\includegraphics[width=0.9\textwidth]{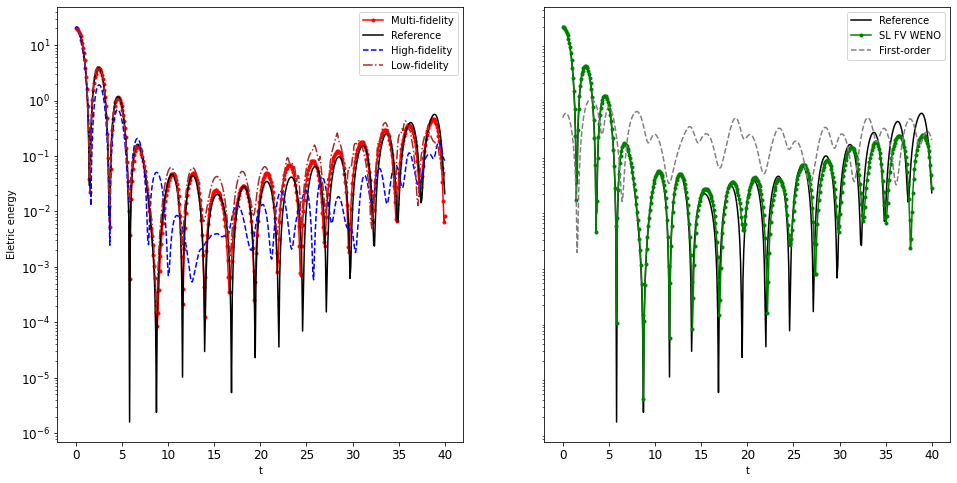}}
	\caption{Time histories of the electric energy of the strong Landau damping in Example \ref{landau-damping} with $\alpha=0.5$.}\label{fig:landau_eletric}
\end{figure}

In Figure \ref{fig:landau_eletric}, we plot the time histories of the electric energy for each approach. Our method with the second RKEI yields results that agree well with the reference solution, significantly outperforming all other methods considered in the testing. 



\begin{example}\label{two-stream}
    In this example, we simulate the symmetric two stream instability with the initial condition 
    \begin{equation}\label{eq:two-stream}
        f(x,v,t=0)=\frac{1}{\sqrt{2\pi}}(1+\alpha \cos(kx))v^2\exp\left(-\frac{v^2}{2}\right),\quad  x\in[0,L],\quad v\in[-V_c,V_c],
    \end{equation}
    where $k=0.5$, $L=4\pi$, and $V_c=2\pi$. 
\end{example}

We generate the low-fidelity training data by obtaining five solution trajectories over a $32\times 64$ grid. To generate the high-fidelity training data, we coarsen two high-resolution solution trajectories over a $256\times 512$-cell grid by a factor of 8 in each dimension. The initial conditions are determined using \eqref{eq:two-stream}, with $\alpha$ randomly sampled from a uniform distribution in the range $[0.01, 0.05]$. We use a CFL number of $1.8$ and set the stencil size to be $5\times 5$.  We choose the loss weights as $\lambda_1=0.2$ and $\lambda_2=1$. The reference solution is generated with the same approach used to create the high-fidelity data. 

\begin{figure}[!htbp]
    \centerline{\includegraphics[width=0.9\textwidth]{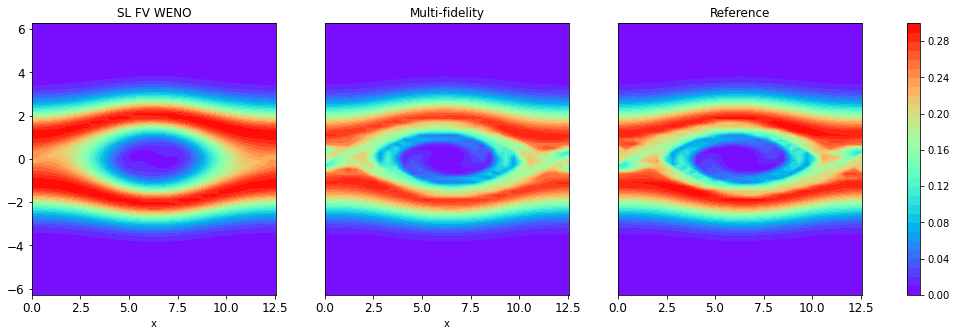}}
	\caption{Contour plots of the numerical solutions for the two stream instability at $t=53$ in Example \ref{two-stream} with $\alpha=0.01$.}\label{fig:two_contour}
\end{figure}

\begin{figure}[!htbp]
    \centerline{\includegraphics[width=0.9\textwidth]{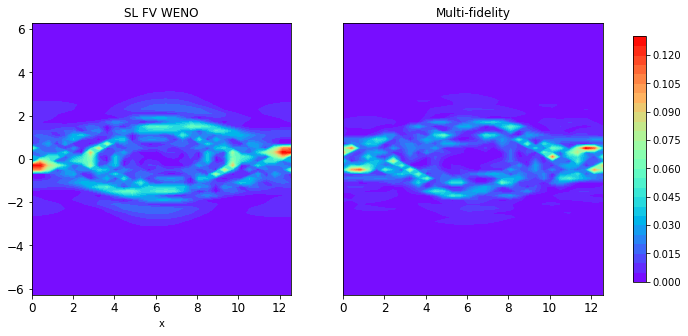}}
	\caption{The error of the numerical solutions for the two stream instability at $t=53$ in Example\ref{two-stream} with $\alpha=0.01$.}\label{fig:two_error}
\end{figure}

For testing, we report contour plots of the two-stream instability with $\alpha=0.01$ at $T=53$ in Figure \ref{fig:two_contour} for the multi-fidelity method as well as the SL FV WENO scheme for comparison.  The proposed multi-fidelity method produces high-quality numerical results that agree well with the reference solution, similar to the previous example. However, the SL FV WENO scheme fails to capture fine-scale structures of interest, such as the roll-up at the center of the solution. We also present the absolute error between the numerical solutions and  the reference solution in Figure \ref{fig:two_error}.



\begin{example}\label{two-stream-another}
In this example,  we consider another two stream instability with the following initial condition
    \begin{equation}
\begin{aligned}\label{eq:two-stream-another} f(x, v, t=0) 
=  \frac{2}{7 \sqrt{2 \pi}}(1+5 v^2)(1+\alpha_1\cos(kx)+\alpha_2\cos(2kx)+\alpha_3\cos(3kx)) \exp \left(-\frac{v^2}{2}\right),
\end{aligned}
\end{equation}
    where $k=0.5$, $x\in[0,4\pi]$, and $v\in[-2\pi,2\pi]$. 
\end{example}

This example is more challenging than the previous two examples, as the initial condition is defined by perturbing  the first three Fourier modes of a equilibrium with magnitudes $\alpha_1$, $\alpha_2$, and $\alpha_3$, respectively.
We generate the low-fidelity training data by obtaining 6 solution trajectories on a $32\times 64$ grid. For the high-fidelity data, we downsample two high-resolution trajectories on a $256\times 512$-cell grid by a factor of 8 in each dimension. The initial conditions for the training data sets are determined using \eqref{eq:two-stream-another} with $\alpha_1, \alpha_2$, and $\alpha_3$ randomly selected from a uniform distribution in the range $[0.01, 0.02]$. For comparison, the reference solution is generated with the same approach used to create the high-fidelity data. We use a CFL number of 1.8 and a stencil size of $5\times 5$.  $\lambda_1$ and $\lambda_2$ in \eqref{eq:multi-loss} are chosen as 0.1 and 1, respectively.

\begin{figure}[!htbp]
    \centerline{\includegraphics[width=0.9\textwidth]{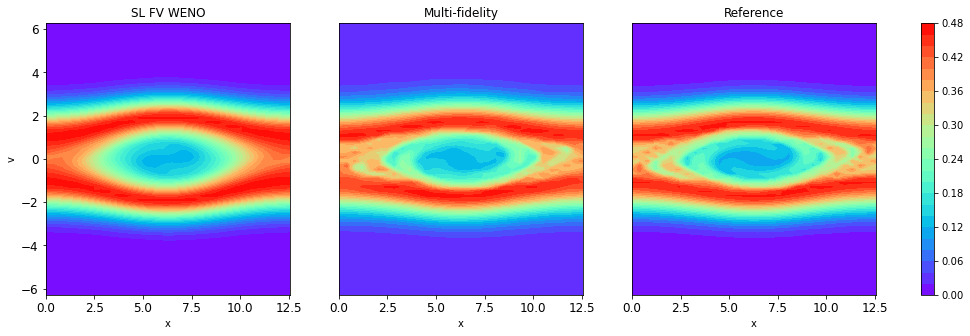}}
	\caption{Contour plots of the numerical solutions of the two stream instability at $t=53$ in Example \ref{two-stream-another} with $\alpha_1=0.01, \alpha_2=0.01/1.2, \alpha_3=0.01/1.2$.}\label{fig:two_contour_another}
\end{figure}

\begin{figure}[!htbp]
    \centerline{\includegraphics[width=0.9\textwidth]{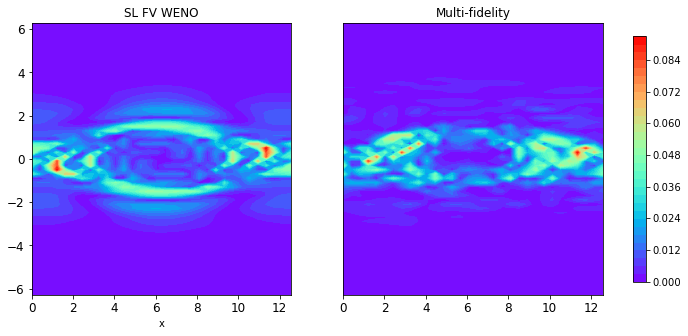}}
	\caption{The error of the numerical solutions of the two stream instability at $t=53$ in Example \ref{two-stream-another} with $\alpha=0.01$.}\label{fig:two_another__error}
\end{figure}

During testing, we consider the initial condition  with  $\alpha_1=0.01, \alpha_2=0.01/1.2, \alpha_3=0.01/1.2$, which is a widely used benchmark configuration in the literature. Note that such a parameter choice is outside the range of the training data.
In Figure \ref{fig:two_contour_another}, we present contour plots of the numerical solutions computed by our multi-fidelity method and the SL FV WENO scheme. It is observed that the results by our method qualitatively agree with the reference solution, effectively capturing the underlying fine-scale structures of interest.
The SL FV WENO scheme, on the other hand, can provide reasonable results, but tends to smear out the small-scale structures  in the solution.  This demonstrates the proposed multi-fidelity model is capable of producing results with reasonable accuracy, highlighting its generalization capabilities.
 Figure \ref{fig:two_another__error} presents the errors between the numerical solutions and the reference solution.
 

\section{Conclusion}\label{sec:con}
In this paper, we have proposed a novel multi-fidelity ML-based  SL FV scheme for solving transport equations. This method is specifically designed for scenarios  where there is an abundance of low-fidelity data and a limited amount of high-fidelity data. By using a composite NN architecture, our method can effectively approximate the inherent correlation between the high-fidelity and low-fidelity data. Numerical experiments  conducted in this study  show that the proposed method achieves improved stability and accuracy 
compared to networks trained solely on either low-fidelity data or high-fidelity data. This indicates that the multi-fidelity approach enhances the performance and capabilities of the SL FV scheme for solving transport equations.
Furthermore, we have extended this multi-fidelity method to the simulation of nonlinear Vlasov-Poisson systems, coupled with the high-order RKEI. This extension allows for accurate and efficient simulations of complex physical phenomena for this multi-fidelity approach. Future work includes exploring the possibility of applying the method to more complicated systems, further improving the generalization capabilities of the method, investigating the use of graph NNs for accommodating unstructured meshes, and address challenges related to adaptivity and complex geometries, among other potential research directions.


\section*{Acknowledgments}
Research work of W. Guo is partially supported by the NSF grant NSF-DMS-2111383, Air Force Office of Scientific Research FA9550-18-1-0257. Research work of X. Zhong is partially supported by the NSFC Grant 12272347.

 \bibliographystyle{abbrv}
 \bibliography{main}
\end{document}